# Space-state complexity of Korean chess and Chinese chess

Donghwi Park
Korea University, Seoul

**Abstract**
This article describes how to calculate exact space-state complexities of Korean chess and Chinese chess. The state-space complexity (a.k.a. search-space complexity) of a game is defined as the number of legal game positions reachable from the initial position of the game.[1][3] The number of exact space-state complexities are not known for most of games. However, we calculated actual space-state complexities of Korean chess and Chinese chess.

## I. Introduction

Chinese chess (Xiangqi) is a popular absolute-strategy board game in China. Korean chess (Janggi) is a popular absolute-strategy board game in Korea, which is derived from Chinese chess. They share very similar pieces and boards. The state-space complexity of a game is defined as the number of legal game positions reachable from the initial position of the game.[1] The number of exact space-state complexities are not known for most of games which include Western chess or Japanese chess (Shogi).[4] There are a few exceptions when number of exact space-state complexities are small, such as Dobutsu shogi[5] or Tic-tac-toe[6]. Unlike Western chess or Japanese chess, Korean chess and Chinese chess do not have promotion rule, which enable calculation of the exact state-space complexity possible.

## II. Chinese chess

Chinese chess board has ten ranks and nine files. There is a river between central two ranks. River splits the board. Each side of has a Palace. Each Palace is located in central three files and outmost three ranks. There are intersections which link corners of each Palace to central point of the palace.[2]

Each side has sixteen pieces: one King(General), two Advisors, two Elephants, two Horses, two Rooks, two Cannons, and five Pawns(Soldiers). Kings move one space horizontally or vertically at a time, not diagonally along the diagonal lines in palaces. Advisors move one space diagonally along the diagonal lines in palaces, meaning they cannot leave their palaces. An elephant moves and captures one point diagonally. They cannot jump over other piece. A horse moves and captures one point horizontally or vertically and then one point diagonally which are not adjacent its former position. A rook moves any distance horizontally or vertically, not diagonally along the diagonal lines in palaces. Cannons move like rooks if they do not capture other pieces. but can only capture by jumping a single piece, along the horizontal or vertical line. Each side has five pawns. a Pawn can only move forward before it cross the river. Pawns have crossed the river, they may also move one point horizontally.[2]

Elephant is restricted to seven positions, Advisor is restricted to five positions, King is restricted to nine positions. Pawns in one side are restricted to fifty-five positions. However pawns cannot move horizontally before it cross the river, every pawn in oneside which have not crossed the river belongs to original files.

**II.1. Arrangement of King, Advisors, and Elephants.**
In Chinese chess, there are no sites occupied by one player's Elephants that can be occupied by the same player's Advisors. Every site that can be occupied by one player's advisor can be occupied by the same player's King. There is one site that can be occupied by one player's King and the same player's Elephants. There are two sites that can be occupied by one player's Pawns and the same player's Elephants. These sites are on the fifth rank and on the third or seventh plies. We will distinguish the count number of the Arrangement by the position of the Elephants. Elephants on fifth rank or third rank and fifth file should be distinguished from other cases.

If we arrange two Advisors on five sites, there are ten cases. If we arrange one Advisor on five sites, there are five cases. If we arrange zero Advisor on five sites, there are one cases.

If we arrange two Elephants on seven sites and they are on the fifth rank, no other Elephants of that player exist. Therefore, there are no Elephants on the palace. There is one case. If we arrange two Elephants on seven sites and one of them is on the fifth rank, the other Elephants of that player exist. There are two cases; the other Elephant is on the third rank and the fifth file, and there are eight cases for elsewhere. If we arrange two Elephants on seven sites and none are on the fifth rank, there are eight cases: another Elephant is on the palace, and there are six cases for elsewhere.

If we arrange one Elephant on seven sites, we can place the one or zero elephant on fifth rank. If we arrange one Elephant on seven sites and that is on the fifth rank, no other Elephant of that player exist. So, there are no Elephant on the palace. There are two cases. If we arrange one Elephant on seven sites and that is not on the fifth rank, there are one case when the Elephant is on the palace and four case when the Elephant is not on third rank&fifth file.

If we arrange zero Elephant on seven sites, there are only one case.

The King should be placed on sites in the palace that are not occupied by Elephants or Advisors.
In one's palace, Elephants or Advisors can occupy zero to three sites.

If we arrange two Elephants and two Advisors, the number of cases is ten times the number of arrangements of two Elephants on seven sites.

If we arrange two Elephants on the fifth rank and two Advisors, no other Elephants of that player exist. Therefore, there are ten cases. If we arrange Elephants on seven sites and one of them is on the fifth rank and arrange two Advisors, there are twenty cases in which the other Elephant is in the palace and eighty cases in which the other Elephant is not in the palace. If we arrange Elephants on seven sites and none are on the fifth rank and arrange two Advisors, there are forty cases in which the other Elephant is in the palace and sixty cases in which the other Elephant is not in the palace.

In these situations, if there is one Elephant in the palace, we can place the King on six sites. if there is no Elephant in the palace, we can place the King on seven sites. If we arrange two Elephants on the fifth rank, we can place the King on seven sites. In this case, the number of arrangements of one's Advisors, Elephants and King on one's camp is 70.
If we place one Elephant on the fifth rank, we can place the King on six sites for 20 cases and on seven sites for 80 cases. In this case, the number of arrangements of one's Advisors, Elephants and King on one's camp is 680.
If we place no Elephant on the fifth rank, we can place the King on six sites for 40 cases and on seven sites for 60 cases. In this case, the number of arrangements of one's Advisors, Elephants and King on one's camp is 660.

If we arrange one Elephant and two Advisors, the number of cases is ten times the number of arrangements of one Elephant on seven sites.
In this case, we cannot arrange two Elephants on the fifth rank. There are zero cases.
If we arrange the Elephant on seven sites and that is on the fifth rank and arrange two Advisors, there are twenty cases in which this Elephant are not in the palace.
If we arrange the elephant on seven sites and that is not on the fifth rank and arrange two Advisors, there are ten cases in which this Elephant is in the palace and forty cases in which this Elephant is not in the palace.

In these situations, if there is one Elephant in the palace, we can place the King on six sites. if there is no Elephant in the palace, we can place the King on seven sites. In this case, we cannot arrange two Elephants on the fifth rank. If we arrange one Elephant on the fifth rank, we can place the King on seven sites. In this case, the number of arrangements of one's Advisors, Elephant and King on one's camp is 140. If we arrange no Elephant on the fifth rank, we can place the King on six sites for 10 cases and on seven sites for 40 cases. In this case, the number of arrangements of one's Advisors, Elephant and King on one's camp is 340.
If we arrange no Elephants and two Advisors, there are no elephant on the fifth rank or in the palace. There are ten cases. We can place the King on seven sites. In this case, the number of arrangements of one's Advisors and King on one's camp is 70.

If we arrange two Elephants and one Advisor, the number of cases is five times the number of arrangements of two Elephants on seven sites.

If we arrange two Elephants on the fifth rank and one Advisor, no other Elephants of that player exist. Therefore, there are five cases. If we arrange Elephants on seven sites and one of them is on the fifth rank and arrange one Advisor, there are ten cases in which the other Elephant is in the palace and forty cases in which the other Elephant is not in the palace. If we arrange Elephants on seven sites and none are on the fifth rank and arrange one Advisor, there are twenty cases in which the other Elephant is in the palace and thirty cases in which the other Elephant is not in the palace.

In these situations, if there is one Elephant in the palace, we can place the King on seven sites. if there is no Elephant in the palace, we can place the King in eight sites. If we arrange two Elephants on the fifth rank, we can place the King on eight sites. In this case, the number of arrangements of one's Advisor, Elephants and King on one's camp is 40.
If we place one Elephant on the fifth rank, we can place the King on seven sites for 10 cases and on eight sites for 40 cases. In this case, the number of arrangements of one's Advisor, Elephants and King on one's camp is 390.
If we place no Elephant on the fifth rank, we can place the King on seven sites for 20 cases and on eight sites for 30 cases. In this case, the number of arrangements of one's Advisor, Elephants and King on one's camp is 360.

If we arrange one Elephant and one Advisor, the number of cases is five times the number of arrangements of one Elephant on seven sites.
In this case, we cannot arrange two Elephants on the fifth rank. There are zero case.
If we arrange Elephant on seven sites and that is on the fifth rank and arrange one Advisor, there are ten cases in which this Elephant are not in the palace. If we arrange the elephant on seven sites and that is not on the fifth rank and arrange one Advisor, there are five cases in which this Elephant is in the palace and twenty cases in which this Elephant is not in the palace.

In these situations, if there is one Elephant in the palace, we can place the King on seven sites. if there is no Elephant in the palace, we can place the King on eight sites. In this case, we cannot arrange any Elephant on the fifth rank. If we arrange one Elephant on the fifth rank, we can place the King on eight sites. In this case, the number of arrangements of one's Advisor, Elephant and King on one's camp is 80. If we arrange no Elephant on the fifth rank, we can place the King on seven sites for 5 cases and on eight sites for 20 cases. In this case, the number of arrangements of one's Advisor, Elephant and King on one's camp is 195.

If we arrange no Elephants and one Advisor, there are no elephant on the fifth rank or in the palace. There are five cases. We can place the King on eight sites. In this case, the number of arrangements of one's Advisor and King on one's camp is 40.

Suppose we arrange two Elephants and no Advisor. If we arrange two Elephants on the fifth rank, no other Elephants of that player exist. Therefore, there are one cases. If we arrange Elephants on seven sites and one of them is on the fifth rank, there are two cases in which the other Elephant are in the palace and eight cases in which the other Elephant are not in the palace. If we arrange Elephants on seven sites and none are is on the fifth rank, there are for cases in which the other Elephant are in the palace and six cases in which the other Elephant are not in the palace.

In these situations, if there is one Elephant in the palace, we can place the King on eight sites. if there is no Elephant in the palace, we can place the King in nine sites. If we arrange two Elephants on the fifth rank, we can place the King in nine sites. In this case, the number of arrangements of one's Elephants and King on one's camp is 9.
If we arrange one Elephant on the fifth rank, we can place the King on eight sites for 2 cases and on nine sites for 8 cases. In this case, the number of arrangements of one's Elephants and King on one's camp is 88.

If we arrange no Elephant on the fifth rank, we can place the King on eight sites for 4 cases and on nine sites for 6 cases. In this case, the number of arrangements of one's Elephants and King on one's camp is 86.

Suppose we arrange one Elephant and no Advisor.

In this case, we cannot arrange two Elephants on the fifth rank. There are zero case.
If we arrange Elephant on seven sites and that is on the fifth rank, there are two cases in which this Elephant are not in the palace. If we arrange the elephant on seven sites and that is not on the fifth rank, there are one cases in which this Elephant is in the palace and four cases in which this Elephant is not in the palace.

In these situations, if there is one Elephant in the palace, we can place the King on eight sites. if there is no Elephant in the palace, we can place the King on nine sites. In this case, we cannot arrange two Elephants on the fifth rank. If we arrange one Elephant on the fifth rank, we can place the King on nine sites. In this case, the number of arrangements of one's Elephant and King on one's camp is 18.
If we arrange no Elephant on the fifth rank, we can place the King on eight sites for 1 cases and on nine sites for 4 cases. In this case, the number of arrangements of one's Elephant and King on one's camp is 44.

If we arrange no Elephants and no Advisors, we can place the King on nine sites. In this case, the number of arrangements of one's King on one's camp is 9.

These tables are comprised of the above results.

| used pieces | Adv/Ele | total count | two 5th Elp | one 5th Elp | no 5th Elp |
|---|---|---|---|---|---|
| 5 | 2,2 | 1410 | 70 | 680 | 660 |
| 4 | 2,1 | 480 | 0 | 140 | 340 |
| 3 | 2,0 | 70 | 0 | 0 | 70 |
| 4 | 1,2 | 810 | 40 | 390 | 380 |
| 3 | 1,1 | 275 | 0 | 80 | 195 |
| 2 | 1,0 | 40 | 0 | 0 | 40 |
| 3 | 0,2 | 183 | 9 | 88 | 86 |
| 2 | 0,1 | 62 | 0 | 18 | 44 |
| 1 | 0,0 | 9 | 0 | 0 | 9 |

Table 1

| used pieces | total count | two 5th Elp | one 5th Elp | no 5th Elp |
|---|---|---|---|---|
| 5 | 1410 | 70 | 680 | 660 |
| 4 | 1290 | 40 | 530 | 720 |
| 3 | 528 | 9 | 168 | 351 |
| 2 | 102 | 0 | 18 | 84 |
| 1 | 9 | 0 | 0 | 9 |

Table 2

### II.2. Arrangement of King, Advisors, Elephants and Soldiers.

Then, let us count the number of arrangements of soldiers on my side. In chinese chess, there are ten sites of one player's soldiers on his side; fourth or fifth rank in odd number files. Soldiers cannot move backward and they cannot move one point horizontally before crossing the river. Therefore, two soldiers cannot exist in a file on my side.

If we arrange my soldiers on my side, we choose files and choose the rank for each soldier. However, in some cases, one or two sites in these ten sites might be occupied by my elephants. This is the result of arranging my soldiers on my side by the number of unoccupied sites in these ten sites and the number of my soldiers on my side.

| soldier \| blank sites | 10 | 9 | 8 |
|---|---|---|---|
| 0 | 1 | 1 | 1 |
| 1 | 10 | 9 | 8 |
| 2 | 40 | 32 | 25 |
| 3 | 80 | 56 | 38 |
| 4 | 80 | 48 | 28 |
| 5 | 32 | 16 | 8 |

Table 3

Additionally, let us count the number of arrangements if we arrange my Advisors, my Elephants, the King, and my soldiers on my side. There are 45 sites on my side.

This table is the number of arrangements of Advisors, my Elephants, the King, and my soldiers on my side. In this table, blank sites (file name) mean sites on my side that are not occupied by Advisors, my Elephants, the King, and my soldiers on my side. Ranks mean the number of my soldiers on my side.

| soldier\|blank | 35 | 36 | 37 | 38 | 39 | 40 | 41 | 42 | 43 | 44 |
|---|---|---|---|---|---|---|---|---|---|---|
| 0 | | | | | | 1410 | 1290 | 528 | 102 | 9 |
| 1 | | | | | 13280 | 12290 | 5094 | 1002 | 90 | |

| | | | | | | | |
|---|---|---|---|---|---|---|---|
| 2 | | | | 49910 | 46760 | 19641 | 3936 | 360 |
| 3 | | | 93540 | 88800 | 37830 | 7728 | 720 |
| 4 | | 87400 | 84160 | 36396 | 7584 | 720 |
| 5 | 32560 | 31840 | 13992 | 2976 | 288 |

Table 4

This new table is the number of arrangements of Advisors, my Elephants, the King, and my soldiers on my side. In this table, blank sites (file name) mean sites on my side that are not occupied by Advisors, my Elephants, the King, and my soldiers on my side. Ranks mean we can arrange same or more soldiers than table's file numbers on the enemy's side.

| soldier\|blank | 35 | 36 | 37 | 38 | 39 | 40 | 41 | 42 | 43 | 44 |
|---|---|---|---|---|---|---|---|---|---|---|
| 5 | | | | | | 1410 | 1290 | 528 | 102 | 9 |
| 4 | | | | | 13280 | 13700 | 6384 | 1530 | 192 | 9 |
| 3 | | | | 49910 | 60040 | 33341 | 10320 | 1890 | 192 | 9 |
| 2 | | | 93540 | 138710 | 97870 | 41069 | 11040 | 1890 | 192 | 9 |
| 1 | | 87400 | 177700 | 175106 | 105454 | 41789 | 11040 | 1890 | 192 | 9 |
| 0 | 32560 | 119240 | 191692 | 178082 | 105742 | 41789 | 11040 | 1890 | 192 | 9 |

Table 5

In this table, let us name the rank the name n, the file name k and the value of the cell d(n,k). d(n,k) is the number of arrangements of my advisors, elephants, King and soldiers on my side, which satisfy there are exactly n unoccupied sites and at least k unused soldiers.

Then let us calculate the number of arrangement if we arrange Advisors, Elephants, Kings, and Soldiers on the board. There are 90 sites.
let n1 be the number of sites unoccupied by A's pieces on A's side, k1 be the number of A's soldiers on B's side, n2 as the number of sites unoccupied by B's pieces on B's side, and k2 be the number of B's soldiers on A's side.

Then let us calculate the number of arrangements if we arrange Advisors, Elephants, Kings, and Soldiers in the board for given n1,k1,n2,k2.

First, A arranges the total n1 number of Advisors, Elephants, Kings, and soldiers on A's side and leaves at least k1 of A's soldiers. The number of arrangements is d(n1,k1). Second, B arranges total n2 number of Advisors, Elephants, Kings, and soldiers on B's side and leaves at least k2 of B's soldiers. The number of ways is d(n2,k2). Third, B arranges B's soldiers to n1 blank sites on A's side. The number of arrangements is c(n1,k2). Forth, A arranges A's soldiers to n2 blank sites on B's side. The number of arrangements is c(n2,k1).

The number of total arrangements is d(n1,k1)×d(n2,k2)×c(n1,k2)×c(n2,k1). In this case, there are n1+n2-k1-k2 sites are unoccupied by any pieces. The boundary condition is 5≥k1,k2≥0, 44≥n1,n2, n1-k1≥35, n2-k2≥35. let n=n1+n2-k1-k2.

The number of arrangements when we arrange the advisors, elephants, soldiers, and Kings with n blank sites is

$$\sum_{\substack{35 \leq n_1, n_2 \leq 44 \\ 0 \leq k_1, k_2 \leq 5 \\ n_1 - k_1 \geq 35 \\ n_2 - k_2 \geq 35 \\ n_1 + n_2 - k_1 - k_2 = n}} d(n_1, k_1) d(n_2, k_2) \binom{n_1}{k_2} \binom{n_2}{k_1}$$

6072015837104228000 (n=70)
13932273683634608000
15302416298575447500
10415675878701420000
4850335101880323628
1620169838558710348
398556758971233856
73409438301988732
10306140239862692
1131080570393880
100447213926298
7330142404440
444595549080
22199620332
900695862
28838016
655542
10584
81 (n=88).
Let us this number as k(n).

## II.3. Arrangement of other pieces and space-state complexity of Chinese chess.

Now, let us calculate the number of arrangements of Chariots, Horses, and Cannons on unoccupied sites. A and B have two chariots, horses and cannons. Therefore, there are six labeled pieces. (See OEIS A141765). The number of distinct ways in which the six labeled pieces can be distribute on n sites allowing at most two pieces to fall on each sites are T(6,n). Let d(n) be T(6,n).

d(0) to d(12) is
1,6,36,210,1170,6120,29520,128520,491400,1587600,4082400,7484400,7484400이다.

The number of distinct ways for arranging pieces on the board when there are x sites unoccupied by any pieces except chariots, horses and cannons and arranging chariots, horse and cannons on to y sites is $k(x)\binom{x}{y}d(y)$. So, space-state complexity of Chinese chess is $\sum_{\substack{70 \leq x \leq 88 \\ 0 \leq y \leq 12}} k(x)\binom{x}{y}d(y)$. The result is
75879095159780903710155382525117211 50667.

## III. Korean chess
Korean chess uses a similar board than Chinese chess. In Korean chess, every piece except Kings, advisors and soldiers can occupy any sites. Soldiers cannot move backward. Therefore, Soldiers can occupy any sites except the first, second and third rank. Sites in the backward three ranks cannot be occupied by A's soldiers, both player's soldiers can occupy the middle four ranks. The player's soldiers cannot occupy Palace of the player.

We will divide the board into three areas. A's home zone; three ranks which include A's palace, middle four ranks, B's home zone; three ranks which include B's palace.

### III.2. Arrangement of King, Advisors and Soldiers in Korean chess.
In one's palace, the King and advisors can occupy any sites in the palace. In these situations, the number of arrangements are 9=c(9,1) for no advisor, 72=2×c(9,2) for one advisor, and 252=3×c(9,3) for two advisors.

Let d(n,k) be the number of arrangements for totally distributing the n pieces of my King, my advisors, and the opponent's soldier on my home zone and at least k opponent's soldiers are unused. Then, $d(n,k) = \sum_{\substack{i=1,2,3 \\ 0 \leq n-i \leq 5 \\ 5-n+i \geq k}} i \times \binom{9}{i} \times \binom{27-i}{n-i}$. In this formula, i is the number of my King or advisors in my palace. Therefore i is 1, 2, 3. The number of arrangements of the King and advisors in someone's palace is $i \times \binom{9}{i}$. Then, n-i opponent's soldiers must be arranged in 27-i unused sites in my home zone. The boundary condition is 5-n+i≧k.

| n\k | 0 | 1 | 2 | 3 | 4 | 5 |
|---|---|---|---|---|---|---|
| 1 | 9 | 9 | 9 | 9 | 9 | 9 |
| 2 | 306 | 306 | 306 | 306 | 306 | 72 |
| 3 | 4977 | 4977 | 4977 | 4977 | 2052 | 252 |
| 4 | 51048 | 51048 | 51048 | 27648 | 6048 | |
| 5 | 369702 | 369702 | 235152 | 69552 | | |
| 6 | 2012868 | 1420848 | 510048 | | | |
| 7 | 6503112 | 2677752 | | | | |
| 8 | 10711008 | | | | | |

Table 6. value of d(n,k)

Then let us calculate the number of arrangements if we arrange Advisors, Kings, and soldiers in the board.
Let n1 be the number of sites occupied by A's advisors , A's King or B's soldiers in A's home zone, n2 be the number of sites occupied by B's advisors, B's King or A's soldiers in B's home zone, k1 be the number of B's soldiers in the four middle ranks, and k2 be the number of A's soldiers in the four middle ranks,

Then, let us calculate the number of arrangements if we arrange Advisors, Kings, and soldiers on the board for given n1, k1, n2, and k2.

First, A arranges the total n1 numbers of A's advisors,  A's King, and B's soldiers on A's home zone and leaves at least k1 of B's soldiers for elsewhere. The number of arrangements is d(n1,k1). Second, B arranges the total n2 numbers of B's advisors,  B's King, and A's soldiers on B's home zone and leaves at least k2 of A's soldiers for elswhere. The number of arrangements is d(n2,k2). Third, B arranges B's k1 soldiers to the four middle ranks, and A arranges A's  k2 soldiers to the four middle ranks, The number of arrangements is c(36,k1)×c(36-k1,k2). The number of these arrangements is d(n1,k1)×d(n2,k2)×c(36,k1)×c(36-k1,k2). In this case, we used n1+k1+n2+k2 pieces. and 1≦n1,n2≦8, 0≦k1,k2≦5.

Let k(n) as the number of arrangements n pieces which are advisors, Kings or soldiers on Janggi board,

$$k(n) = \sum_{\substack{1 \le n_1, n_2 \le 8 \\ 0 \le k_1, k_2 \le 5 \\ n_1 + n_2 + k_1 + k_2 = n}} d(n_1, k_1) d(n_2, k_2) \binom{36}{k_1} \binom{36 - k_1}{k_2}.$$

Let us calculate this formula. From k(16) to k(2).
1457601002568716544,1185971655381537024,470042212117883328,111504273140075328,17627589996960672,1967816967471936,171617399962470,12043618055460,686813883426,31907861496,1200808557,35663652,783918,11340,81.

### III.2. Arrangement of other pieces and Space-state complexity of Korean chess

Now, let us calculate the number of arrangements of Chariots, Horses, Elephants and Cannons on unoccupied sites. A and B have two chariots, horses, elephants and cannons. Therefore, there are eight labeled pieces. (See OEIS A141765). The number of distinct ways in which the eight labeled pieces can be distributed on n sites allowing at most two pieces to fall in each sites are T(8,k). Let s(k) be T(8,k). s(0)=1, s(1)=8. $s(k) = \sum_{\substack{k-i \leq 8 \\ k-2i \geq 0 \\ i \geq 0}} \binom{8}{i} \times \binom{8-i}{k-2i} \times \frac{k!}{2^n}$ and this results are

1,8,64,504,2028,28560,44520,294000,441840,6773760,6827940,209933640,209766060,5448713760,5448660840,40864824000,40864824000 (from s(0) to s(16) )

The number of whole possible arrangements in Korean chess is $\sum_{n=2}^{16} \sum_{k=0}^{16} s(k) k(n)$.

The final result(State-space complexity of Korean chess) is
23510395465980130401868412314878554298901846.

### IV. Reference

[1]Victor Allis, Searching for Solutions in Games and Artificial Intelligence, Ph.D. thesis, University of Limburg, Maastricht, The Netherlands, 1994. pp. 158
[2]Shi-Jim Yen, Jr-Chang Chen, Tai-Ning Yang, and Shun-Chin Hsu. Computer Chinese Chess, International Computer Games Association Journal 27 (1), 2004. pp. 3-18
[3]H. Jaap van den Herik, Jos W.H.M. Uiterwijk, Jack van Rijswijck.
Games solved: Now and in the future, Artificial Intelligence 134, 2002, pp. 277-311
[4]Hitoshi Matsubara and Reijer Grimbergen, Differences between Shogi and western Chess from a computational point of view, Board Games in Academia, an interdisciplinary, Leiden, The Netherlands , 1997
[5]Tetsuro Tanaka, An Analysis of a Board Game "Doubutsu Shogi" (Japanese) IPSJ SIG Notes, Vol. 2009-GI-22, No. 3, 2009. pp. 1-8 (Japanese)
[6]Susan L. Epstein, Learning in the Right Places, Journal of the Learning Sciences 4(3), 1993, pp. 281-319